\documentclass{amsart}

\usepackage{amsmath}
\usepackage{paralist}
\usepackage{graphicx}  
\usepackage{hyperref}

\usepackage[caption=false]{subfig}
\usepackage{enumerate}

\theoremstyle{definition}

\newcommand{\argmin}{\text{argmin}}
\newcommand{\R}{\mathbb{R}}
\newcommand{\TV}{\text{TV}}
\newcommand{\norm}[1]{\left\| #1 \right\|}

\newcommand{\diag}{\text{diag}}
\newcommand{\soft}{\text{soft}}

\newcommand{\dct}{\text{DCT}}
\newcommand{\dctsec}{\text{DCT}^{\text{II}}}
\newcommand{\dctthi}{\text{DCT}^{\text{III}}}
\newcommand{\abs}[1]{\left| #1 \right|}
\newcommand{\sgn}{\text{sgn}}

\title[TV-regularized NMF for smooth unmixing] 
      {Total variation regularized non-negative matrix factorization for smooth hyperspectral unmixing}

\author[Adrien Faivre and Cl\'ement Dombry]{}



\thanks{$^*$ Corresponding author: Adrien Faivre}

\begin{document}
\maketitle

\centerline{\scshape Adrien Faivre$^*$}
\medskip
{\footnotesize
 \centerline{Digital Surf, 16 Rue Lavoisier, 25000 Besan{\c c}on, FRANCE}
	\centerline{Email: afaivre@digitalsurf.fr}
   \centerline{and}
   \centerline{Laboratoire de Math\'ematiques de Besan\c{c}on, UMR CNRS 6623,}
 	\centerline{Universit\'e de Bourgogne Franche-Comt\'e,}
	\centerline{16 route de Gray, 25030 Besan{\c c}on cedex, FRANCE}
} 

\medskip

\centerline{\scshape Cl\'ement Dombry}
\medskip
{\footnotesize
   \centerline{Laboratoire de Math\'ematiques de Besan\c{c}on, UMR CNRS 6623,}
 	\centerline{Universit\'e de Bourgogne Franche-Comt\'e,}
	\centerline{16 route de Gray, 25030 Besan{\c c}on cedex, FRANCE}
	\centerline{Email: clement.dombry@univ-fcomte.fr}
}

\bigskip

\begin{abstract}
Hyperspectral analysis has gained popularity over recent years as a way to infer 
what materials are displayed on a picture whose pixels consist of a mixture of spectral 
signatures. Computing both signatures and mixture coefficients is known 
as unsupervised unmixing, a set of techniques usually based on non-negative matrix 
factorization. Unmixing is a difficult non-convex problem, and algorithms may 
converge to one out of many local minima, which may be far removed from the true global 
minimum. Computing this true minimum is NP-hard and seems therefore out of reach. 
Aiming for interesting local minima, we investigate the addition of total variation 
regularization terms. Advantages of these regularizers are two-fold. 
Their computation is typically rather light, and they are deemed to preserve sharp edges 
in pictures. This paper describes an algorithm for regularized hyperspectral unmixing 
based on the Alternating Direction Method of Multipliers.
\end{abstract}

\section{Hyperspectral unmixing}

Spectral sensors can nowadays acquire electromagnetic intensity with a resolution of 
thousands of wavebands, spread across an increasing number of pixels.  
Analyzing hyperspectral data is therefore of growing complexity. This article tackles 
with hyperspectral unmixing, a matrix factorization technique designed to retrieve 
spectral signatures of pure materials and corresponding proportions from an hyperspectral 
image, while enhancing smoothness of both spectrums and mixing proportions. Pure spectral 
signatures and their mixing proportions are usually respectively referred to as endmembers 
and abundances throughout the hyperspectral analysis literature.

Spectral signatures are by nature positive, and abundances are by definition positive 
proportions. These facts lead to a formulation of unmixing as the following constrained 
optimization problem:
\begin{equation}
\min_{W \in \R_+^{m \times k}, H \in \R_+^{k \times n}} \frac{1}{2}\| M - W H \|_F^2,
\label{eq:nmf}
\end{equation}
designed to recover $k$ endmembers $W_1, \cdots, W_k$ and their corresponding abundances 
$H_1, \cdots, H_n$ from $n$ observed spectrums $M_1, \cdots, M_n$. 
Notations $W_i, H_i$ and $M_i$ respectively stand for the $i$-th column of matrices $W, H$ and $M$.
The non-negative reals are denoted $\R_+$, and $\R_+^{m \times n}$ stands for the set of $m \times n$ 
non-negative matrices. 
The norm $\| \cdot \|_F$ is the standard Frobenius matrix norm.

Problem \eqref{eq:nmf} is known as Non-negative Matrix Factorization (NMF), and is deemed quite difficult.
While optimizing with respect to $W$ or $H$ alone is a simple non-negative least squares 
problem that can be solved in polynomial time, optimizing with both $W$ and $H$ 
simultaneously turns out to be very hard. Vavasis \cite{vavasis2009complexity} proved 
that deciding whether the non-negative rank of a matrix is the same as its rank is 
NP-complete. Arora et al. \cite{arora2012computing} showed how under the Exponential 
Time Hypothesis, there can be no exact algorithm for NMF running in time $(m n)^{o(r)}$. 
There already exists many algorithms producing approximate solutions to NMF \cite{gillis2014and}. 
They range from simple alternating algorithms \cite{lee2001algorithms}, to more theoretically 
involved convexifications \cite{krishnamurthy2012convex}. In hyperspectral literature, 
one recurring trend is to try to recover first the endmembers, and only then compute 
abundances. This is the path followed by VCA \cite{nascimento2005vertex}, SPA \cite{gillis2014and}, 
SIVM \cite{bauckhage2014purely}, and a few others. 
Another approach is to add constraints to the NMF problem, aiming for more
 interesting local minima. The most obvious one we can add is for columns of $H$ to sum to $1$. 
Abundances are indeed supposed to be proportions. Every column of the abundances matrix must 
therefore belong to the $k$-simplex. We note:

\begin{equation*}
H = \begin{pmatrix}H_1 & H_2 & \cdots & H_n\end{pmatrix} \in \Delta_k^n .
\label{eq:simplex}
\end{equation*}

Hoyer \cite{hoyer2004non} remarked that the matrices resulting from NMF are usually sparse, 
and emphasized this by adding sparseness constraints on $W$ or $H$. More recently, it was considered 
natural for abundances to display some sort of spatial regularity across the image. Inspired 
by Total Variation (TV) regularization techniques, Iordache et al. \cite{iordache2012total} penalized 
the $l_1$ norm of the abundances gradient. A related idea from Warren and 
Osher \cite{warren2015hyperspectral} is to consider that endmembers should also display limited total 
variation. The main topic of the present paper is to combine both spectral and spatial TV 
regularization to produce a smoothed NMF.

Adding all these constraints and regularizers together yields a goal functional a little more 
complicated than the one described in Equation \eqref{eq:nmf}. The TV regularized NMF problem reads
\begin{equation}
\min_{W \in \R_+^{m \times k}, H \in \Delta_k^{n}} \frac{1}{2}\| M - W H \|_2^2 + \mu \TV ( H ) + \lambda \TV( W ).
\label{eq:NMFTV}
\end{equation}
This problem can fortunately be split into easier sub-problems, through techniques such as 
the Alternating Direction Method of Multipliers (ADMM). See the 
seminal paper of Boyd et al. \cite{boyd2011distributed} for an in-depth treatment of ADMM 
techniques. 

The ADMM algorithm allows for the optimization of composite problems of the form
\begin{equation}
\min f(x) + g(z) \text{ such that } A x + B z = c .
\label{eq:admm}
\end{equation}
To achieve this, first form the augmented Lagrangian 

\begin{equation*}
L_{\rho}(x,z,y) = f(x)+g(z)+y^T(Ax+Bz-c)+\frac{\rho}{2}\| A x + B z - c \|_2^2 .
\label{eq:al}
\end{equation*}

Then, using updates
\begin{equation}
\begin{cases}
&x_{k+1} = \argmin_x L_{\rho}(x,z_k,y_k) ,\\
&z_{k+1} = \argmin_z L_{\rho}(x_{k+1},z,y_k) ,\\
&y_{k+1} = y_k + \rho (Ax_{k+1}+Bz_{k+1}-c) ,
\end{cases}
\label{eq:admm_upd}
\end{equation}
one is guaranteed to get a local optimum for problem \eqref{eq:admm} under relatively mild assumptions. 
In this paper we will instead use the equivalent shortened version

\begin{equation*}
\begin{cases}
&x_{k+1} = \argmin_x f(x) + \frac{\rho}{2}\| Ax+Bz_k-c+u_k \|_2^2 ,\\
&z_{k+1} = \argmin_z g(z) + \frac{\rho}{2}\| Ax_{k+1}+Bz-c+u_k \|_2^2 ,\\
&u_{k+1} = u_k + (Ax_{k+1}+Bz_{k+1}-c),
\end{cases}
\label{eq:admm_scaled}
\end{equation*}

that follows from Equation \eqref{eq:admm_upd} by introducing $u_k = y / \rho$.

NMF hardly fits the ADMM framework, but Zhang \cite{zhang2010alternating} studied the convergence
of the following splitting scheme
\begin{equation}
\min_{X_1, X_2, Z_1, Z_2} \| M - X_1 Z_2 \|_F^2 + \chi_+(X_2) + \chi_+(Z_1) \text{ such that }
\begin{cases} X_1 = Z_1, \\ X_2 = Z_2. \end{cases}
\label{eq:ADMM_NMF}
\end{equation}
In Equation \eqref{eq:ADMM_NMF}, $\chi_+$ denotes the 
characteristic function of matrices with positive entries, valued $0$ for matrices with no 
negative entries, and $\infty$ for any other. 
The sub-problem in $X_1$ is a least square problem, and the one in $X_2$ can be 
solved using a projection on the positive orthant. Sub-problems in $Z_1, Z_2$ can be dealt with similarly.
These four sub-problems can be solved efficiently. The difference with standard ADMM formulation is that it does not 
seem possible to properly split Equation \eqref{eq:nmf} as the sum of two functions. 

One advantage of the ADMM technique is that adding regularization terms is rather straightforward, 
enabling us to interleave TV regularization steps between every NMF factor update.
ADMM algorithms were already shown to be quite successful for TV denoising problems \cite{rudin1992nonlinear}.

In Section \ref{sec:TV}, we present how classical ADMM algorithms for image TV 
denoising can be adapted to hyperspectral images, using realistic boundary conditions.
In Section \ref{sec:NMF}, we derive an efficient ADMM algorithm for TV regularized NMF.

\section{Total Variation Denoising} \label{sec:TV}

\subsection{Neuman boundary conditions}

The penalization of total variation is a popular technique used to remove noise, introduced 
by Rudin et al \cite{rudin1992nonlinear}. The idea is that a pure signal should have limited 
variation. Typically, the norm of the discrete gradient is penalized, while faithfulness to 
the original signal is promoted.
Corresponding optimization problem for $1$-dimensional signal $y \in \R^n$ reads
\begin{equation}
\min_x \frac{1}{2} \| y - x \|_2^2 + \lambda \| D_n x \|_p,
\label{eq:rof1d}
\end{equation}
where $D_n$ denotes a discrete differentiation operator for signals of length $n$, and $\lambda$ 
is a parameter controlling the strength of the denoising. Using $p=1$ leads to anisotropic denoising, 
and $p=2$ to isotropic denoising. We will focus on the case $p=1$. Problem \eqref{eq:rof1d} can be solved 
efficiently using ADMM. It can indeed be split as

\begin{equation*} 
\min_{x,z} \frac{1}{2} \| y - x \|_2^2 + \lambda \| z \|_1 \text{ such that } D_n x = z .
\end{equation*}

The sub-problem for $z$ reads
\begin{equation}
\min_z \lambda \| z \|_1 + \frac{\rho}{2} \norm{D_n x - z + u}_{2}^2,
\label{eq:z_upd}
\end{equation}
and can be solved efficiently. Let $\soft( x, \lambda) = \sgn(x) (\abs{x}-\lambda)_+$ denote the 
soft-thresholding operator for real $x$. The solution of sub-problem \eqref{eq:z_upd} is then 
given by

\begin{equation*} 
z_{k+1} = \soft \left( D_n x + u, \frac{\lambda}{\rho} \right),
\end{equation*}

where $\soft$ is applied coefficient-wise.
The $x$ update computation is harder, as it requires to solve
\begin{equation}
(I_n + \rho D_n^T D_n) x = y + \rho ( z - u ), 
\label{eq:inv1d}
\end{equation}
which can be quite difficult for large values of $n$.
One common way to handle this is to pretend that the signal has some sort of specific boundaries. 
Periodic boundaries are for instance often used as they allow to represent operator $D_n$ as a circulant matrix. 
This enables the use of Fast Fourier Transform (FFT) to solve Equation \eqref{eq:inv1d},
which is very efficient in practice. 
However, real images are usually not periodic, and we will therefore focus on slightly more 
complicated boundary conditions called Neuman boundaries.

Those conditions express that the gradient is zero on one side of the signal.
Corresponding gradient operator reads

\begin{equation*}
D_n = \begin{pmatrix}
1 & -1 & & \\
& 1 & -1 & & \\
& & \ddots & \ddots \\
& & & 1 & -1 
\end{pmatrix}.
\label{eq:gradient}
\end{equation*}

In order to solve Equation \eqref{eq:inv1d}, we need to invert matrix $L_n = D_n^T D_n$, which is 
a lengthy computation using standard dense matrix inversion algorithms. But $L_n$ is a
tridiagonal matrix with a very specific structure
\begin{equation}
L_n = \begin{pmatrix}
1 & -1 & & & \\
-1 & 2 & -1 & & \\
& \ddots & \ddots & \ddots & \\
& & -1 & 2 & -1 \\
& & & -1 & 1 \\
\end{pmatrix}.
\label{eq:squared_gradient}
\end{equation}

One key property enabling us to inverse $L_n$ swiftly is that Toeplitz-plus-Hankel matrices can be 
diagonalized by Discrete Cosine Transform (DCT) operators, as shown by Ng et al. \cite{ng1999fast}. 
There are many flavors of DCT. We will focus on DCT-II, as it is the most common one. The DCT-II 
matrix is given by

\begin{equation*}
(\dctsec_n)_{i,j} = \cos i (j + \frac{1}{2}) \frac{\pi}{n} ,
\end{equation*}

which is not orthogonal. Its inverse is known as the DCT-III matrix. The diagonalization of the 
squared gradient operator $L_n$ from Equation \eqref{eq:squared_gradient} reads

\begin{equation*}
\dctsec_n L_n \dctthi_n = \diag (s_1, \cdots , s_n),
\end{equation*}

where the diagonal values $s_1, \cdots , s_n$ can be recovered from

\begin{equation*}
\dctsec_n L_n e_1 = \diag (s_1, \cdots , s_n) \dctsec_n e_1,
\end{equation*}

with $e_1$ the first vector of the canonical basis of $\R^n$.
We get
\begin{equation}
s_i = \frac{\cos(\frac{\pi i}{2 n}) - \cos(\frac{3 \pi i}{2 n})}{\cos(\frac{\pi i}{2 n})} .
\label{eq:diagvalues}
\end{equation}

\subsection{Discrete differentiation in \texorpdfstring{$3$}{3}-dimensions}

Hypercubes are multi-indices tables, and a few notions of tensor calculus can be useful to derive 
update rules for TV regularization. See Appendix for a brief overview of those rules.
The proposed unmixing algorithm is based on spatial, or even spectral-spatial regularization.
For one-dimensional problems we introduced operator $D_n$ in Section \ref{sec:TV}.

\begin{equation*}
D_n x = \left(x_i-x_{i+1}\right)_{1 \leq i \leq n}.
\end{equation*}

Higher orders come with more intricate formulas, especially when unfolding the cube. A synthetic 
and efficient way to write down multidimensional discrete differentiation is through tensor contraction.
For notational simplicity, we focus on the $3$-dimensional case that corresponds to hyperspectral 
images. Our discussion remains valid for any other dimension.
For a cube $\mathcal{Y} \in \R^{m \times n \times o}$, differentiation on the first mode can be written as

\begin{equation*}
\mathcal{Y}  \times_1 D_m = D_m Y^{(1)},
\end{equation*}

where $Y^{(1)}$ denotes the first unfolding of $\mathcal{Y}$ (see Appendix for more details).

Matricization enables us to simply state $3$-dimensional anisotopic TV as

\begin{equation*}
\TV (\mathcal{Y}) = \left\| \mathcal{Y} \times_1 D_m \right\|_1 + \left\| \mathcal H \times_2 D_n \right\|_1 + \left\| \mathcal H \times_3 D_o \right\|_1 .
\label{eq:3ddiff}
\end{equation*}

Combining matricization with ADMM theory, we then derive a simple spectral-spatial filtering method close in spirit 
to the one described by Aggarwal \cite{aggarwal2016hyperspectral}. But instead of using a conjugate-gradient type method for solving 
sparse linear equations every other update, we use two $3$-dimensional DCT.

Such a filter could for instance read

\begin{equation*}
\min_{\mathcal X} \frac{1}{2} \left\| \mathcal Y - \mathcal X \right\|_F^2 + \lambda_s \left\| \mathcal X \times_1 D_m \right\|_1 + \lambda_s \left\| \mathcal X \times_2 D_n \right\|_1 + \lambda_t \left\| \mathcal X \times_3 D_o \right\|_1 .
\label{eq:TVfilter}
\end{equation*}

We split last problem as

\begin{equation*}
\min_{\mathcal X, \mathcal Z_1, \mathcal Z_2, \mathcal Z_3} \frac{1}{2} \left\| \mathcal Y - \mathcal X \right\|_F^2 + \lambda_s \left\| \mathcal{Z}_1 \right\|_1 + \lambda_s \left\| \mathcal{Z}_2\right\|_1 + \lambda_t \left\| \mathcal{Z}_3 \right\|_1
\label{eq:TVfilterSplit}
\end{equation*}

subject to constraints

\begin{equation*}
\begin{cases}
\mathcal{Z}_1 = \mathcal X \times_1 D_m, \\
\mathcal{Z}_2 = \mathcal X \times_2 D_n, \\
\mathcal{Z}_3 = \mathcal X \times_3 D_o .
\end{cases}
\label{eq:TVfilterSplitClauses}
\end{equation*}

In order to compute the zeros the gradient with respect to $\mathcal X$, we need to solve a generalized Sylvester 
equation
\begin{equation}
\mathcal X + \rho \mathcal X \times_1 L_m + \rho \mathcal X \times_2 L_n + \rho \mathcal X \times_3 L_o = \mathcal B ,
\label{eq:Neumann3}
\end{equation}
where

\begin{equation*}
\mathcal B = \mathcal Y
+ \rho (\mathcal{Z}_1 - \mathcal{B}_1) \times_1 D_m^T
+ \rho (\mathcal{Z}_2 - \mathcal{B}_2) \times_1 D_n^T
+ \rho (\mathcal{Z}_3 - \mathcal{B}_3) \times_1 D_o^T .
\label{eq:B}
\end{equation*}

The easiest way to solve the classical Sylvester matrix equation is to vectorize the unknown and use 
Kronecker products to produce an equivalent one. Applying this idea to Equation \eqref{eq:Neumann3} yields
\begin{equation}
A \text{ vec}\left(X^{(1)} \right) = \text{vec}\left(B^{(1)} \right) ,
\label{eq:updateX}
\end{equation}
with

\begin{equation*}
A = \left[ I_{mno} + \rho \left(I_o \otimes I_n \otimes L_m + I_o \otimes L_n \otimes I_m + L_o \otimes I_n \otimes I_m \right) \right] .
\label{eq:bigA}
\end{equation*}

The interested reader can report to the Appendix to see details concerning the previous operation.
Matrix $A$ is potentially enormous and solving Equation \eqref{eq:updateX} may seem impossible at first. However, 
noticing that $A$ can be written as a Kronecker sum $I_{mno} + \rho L_o \oplus L_n \oplus L_m$, 
one can diagonalize $A$ as

\begin{equation*}
A = \dctsec_{o,n,m} \left[I_{mno} + \rho (S_o \oplus S_n \oplus S_m ) \right] \dctthi_{o,n,m},
\label{eq:Adiag}
\end{equation*}

where $\dctsec_{o,n,m}$ is the $3$-dimensionnal DCT-II defined as 
$\dctsec_o \otimes \dctsec_n \otimes \dctsec_m$, and $\dctthi_{o,n,m}$ is defined similarly. 
Matrices $S_o, S_n$ and $S_m$ are diagonal matrices whose values are given by Equation \eqref{eq:diagvalues}.

Assuming Neumann boundary conditions, we can now solve efficiently Equation \eqref{eq:Neumann3} using 
two $3$-dimensional DCT only.
\begin{equation}
\mathcal X = \dctthi_{o,n,m} \left( \frac{\dctsec_{o,n,m} \left( \mathcal B \right)}{\rho \mathcal S + 1} \right),
\label{eq:updateX3}
\end{equation}
where $\mathcal S$ is defined as in \eqref{eq:diagvalues} as
\begin{equation}
s_{i,j,k} = \frac{\cos(\frac{\pi i}{2 m}) - \cos(\frac{3 \pi i}{2 m})}{\cos(\frac{\pi i}{2 m})} +
\frac{\cos(\frac{\pi j}{2 n}) - \cos(\frac{3 \pi j}{2 n})}{\cos(\frac{\pi j}{2 n})} +
\frac{\cos(\frac{\pi k}{2 o}) - \cos(\frac{3 \pi k}{2 o})}{\cos(\frac{\pi k}{2 o})}.
\label{eq:DforupdateX3}
\end{equation}

Tensors $\mathcal Z_1, \mathcal Z_2, \mathcal Z_3$ are then readily updated according to

\begin{equation*}
\begin{cases}
\mathcal Z_1 = \soft \left(\mathcal X \times_1 D_m + \mathcal{U}_1, \frac{\lambda_s}{\rho} \right), \\
\mathcal Z_2 = \soft \left(\mathcal X \times_2 D_n + \mathcal{U}_2, \frac{\lambda_s}{\rho} \right), \\
\mathcal Z_3 = \soft \left(\mathcal X \times_3 D_o + \mathcal{U}_3, \frac{\lambda_t}{\rho} \right),
\end{cases}
\label{eq:updateZ3}
\end{equation*}

where the soft-thresholding operator is applied coefficient-wise.
The $\mathcal{X}$ and $\mathcal{Z}$ updates are relatively simple. The hardest one, the $\mathcal{X}$
update, is computable in $O(mno \log(mno))$. Despite consisting of these simple steps, 
the algorithm is rather slow for big cubes. We improve on this issue in Section \ref{sec:NMF}
by regularizing a smaller tensor than $\mathcal{X}$ obtained through factorization.

\section{Non-negative matrix factorization} \label{sec:NMF}

TV regularization is able to preserve discontinuities in images and is often used in image treatment 
applications, for instance when facing debluring problems \cite{yang2009fast, ng1999fast} that would otherwise 
be numerically unstable.
NMF is an ill-posed problem \cite{gillis2014and} for which solutions are usually not unique. We hope 
that adding TV regularization to the problem can lead to better solutions.
In this section, we show that TV denoising can be integrated to the unmixing phase without slowing 
down convergence significantly.

We have to modify slightly what was introduced in Section \ref{sec:TV} in order to fit 
the NMF framework.
Indeed, the total variation on the third mode does not mean anything any meaning any more. It could increase 
or decrease simply by changing the order in which we stack the endmembers in matrix $W$. That is why
we will not focus on the spectral variation of abundances, but will rather add a TV-regularization term
on each endmember instead. Abundances will be regularized slice by slice with $2$-dimensional TV.

More precisely, we are going to devise an algorithm able to solve the following problem:

\begin{equation*}
\min_{\substack{\mathcal H \in \Delta_k^{m \times n} \\ W \in \R_{+}^{o \times k}}} \frac{1}{2} \left\| \mathcal M - \mathcal H \times_3 W \right\|_F^2 + \lambda_s \left\| \mathcal H \times_1 D_m \right\|_1 + \lambda_s \left\| \mathcal H \times_2 D_n \right\|_1 + \lambda_t \left\| D_o W \right\|_1 ,
\label{eq:NMFTV3d}
\end{equation*}

which is a $3$-dimensional equivalent to problem \eqref{eq:NMFTV}.

As in Section \ref{sec:TV}, we start by splitting the problem as:
\begin{equation}
\min \frac{1}{2} \left\| \mathcal M - \mathcal Z_1 \times_3 X_1 \right\|_F^2 + \lambda_s \left\| \mathcal Z_2 \right\|_1 + \lambda_s \left\| \mathcal Z_3 \right\|_1 + \chi_+\left(\mathcal Z_4\right) + \lambda_t \left\|X_2 \right\|_1 + \chi_+\left(X_3\right), 
\label{eq:NMFTVsplit}
\end{equation}
with the following constraints:

\begin{equation*}
\begin{cases} X_1 = Z_0, \\ X_2 = D_o Z_0, \\ X_3 = Z_0 ,\end{cases} \text{\ \ \ \ \ and\ \ \ \ \ }
\begin{cases} \mathcal{Z}_1 = \mathcal X_0 \\ \mathcal{Z}_2 = \mathcal X_0 \times_1 D_m, \\ \mathcal{Z}_3 = \mathcal X_0 \times_2 D_n, \\ \mathcal{Z}_4 = \mathcal X_0 .\end{cases}
\label{eq:NMFTVsplitclauses}
\end{equation*}

We now list every update for solving problem \eqref{eq:NMFTVsplit}.
\begin{itemize}
\item The first sub-problem to solve is

\begin{equation*}
\min_{\mathcal X_0} \begin{pmatrix}
& \frac{\rho}{2} \left\| \mathcal{Z}_1 - \mathcal X_0 + \mathcal U_4 \right\|_F^2 + \\
& \frac{\rho}{2} \left\| \mathcal{Z}_2 - \mathcal X_0 \times_1 D_m + \mathcal U_5 \right\|_F^2 + \\
& \frac{\rho}{2} \left\| \mathcal{Z}_3 - \mathcal X_0 \times_2 D_n + \mathcal U_6 \right\|_F^2 + \\
& \frac{\rho}{2} \left\| \mathcal{Z}_4 - \mathcal X_0 + \mathcal U_7 \right\|_F^2
\end{pmatrix}.
\end{equation*}

Nullity of gradient implies

\begin{equation*}
2 \mathcal X_0 + \mathcal X_0 \times_1 D_m D_m^T + \mathcal X_0 \times_2 D_n D_n^T = \mathcal B,
\label{eq:updateX0}
\end{equation*}

with $\mathcal B = \mathcal{Z}_1 + \mathcal U_4 +
\mathcal{Z}_2 \times_1 D_m^T + \mathcal U_5 +
\mathcal{Z}_3 \times_2 D_n^T + \mathcal U_6 +
\mathcal{Z}_4 + \mathcal U_7 $.
As in Section \ref{sec:TV}, we use a discrete cosine transform to solve this Sylvester equation.
We get
\begin{equation}
\mathcal X_0 = \dctthi_{k,n,m} \left( \frac{\dctsec_{k,n,m} \left( \mathcal B \right)}{\rho \mathcal S + 1} \right),
\label{eq:DCTupdateX0}
\end{equation}
with $s_{i,j,k} = s_{i, j, 0}$ adapted to the $2$-dimensional case from Equation \eqref{eq:DforupdateX3}.

\item The second sub-problem is an order-$3$ least-squares problem, and reads

\begin{equation*}
\min_{X_1} \frac{1}{2} \left\| \mathcal M - \mathcal Z_1 \times_3 X_1 \right\|_F^2 + \frac{\rho}{2} \left\| X_1 - Z_0 +  U_1 \right\|_F^2 .
\end{equation*}

A third mode unfolding of the previous equation gets us a classical least squares problem

\begin{equation*}
\min_{X_1} \frac{1}{2} \left\| M^{(3)} - X_1 Z_1^{(3)} \right\|_F^2 + \frac{\rho}{2} \left\| X_1 - Z_0 + U_1 \right\|_F^2 . 
\end{equation*}

The update is then given by

\begin{equation*}
X_1 = \left( M^{(3)}{Z_1^{(3)}}^T + \rho \left( Z_0 - U_1 \right) \right) \left( Z_1^{(3)} {Z_1^{(3)}}^T + \rho I_k \right)^{-1} .
\label{eq:updateX1}
\end{equation*}

\item The third sub-problem reads:

\begin{equation*} 
\min_{X_2} \lambda_t \left\|X_2 \right\|_1 + \frac{\rho}{2} \left\| X_2 - D_o Z_0 + U_2 \right\|_F^2 .
\end{equation*}

It is a simple proximal problem, with a straightforward solution, that reads

\begin{equation*}
X_2 = \soft \left( D_o Z_0 - U_2, \frac{\lambda_t}{\rho} \right)
\label{eq:updateX2}
\end{equation*}

\item The fourth sub-problem reads

\begin{equation*}
\min_{X_3} \chi_+(X_3) + \frac{\rho}{2} \left\| X_3 - Z_0 + U_3 \right\|_F^2.
\end{equation*}

The solution is given by the projection $\Pi_+$ on the positive orthant applied coefficient-wise
to the difference between $Z_0$ and $U_3$:

\begin{equation*}
X_3 = \Pi_+ \left( Z_0 - U_3 \right).
\label{eq:updateX32}
\end{equation*}

\item The fifth sub-problem reads

\begin{equation*}
\min_{Z_0} \frac{\rho}{2} \left\| X_1 - Z_0 + U_1 \right\|_F^2 +
\frac{\rho}{2} \left\| X_2 - D_o Z_0 + U_2 \right\|_F^2 +
\frac{\rho}{2} \left\| X_3 - Z_0 + U_3 \right\|_F^2 .
\end{equation*}

It is solved by

\begin{equation*}
Z_0 = \left(D_o^T D_o + 2 I_k\right)^{-1} \left(X_1 + U_1 + D_o^T (X_2 + U_2) + X_3 + U_3 \right).
\label{eq:updateZ0}
\end{equation*}

\item The sixth sub-problem is again a least squares-problem :

\begin{equation*}
\min_{\mathcal Z_1} \frac{1}{2} \left\| \mathcal M - \mathcal Z_1 \times_3 X_1 \right\|_F^2 + \frac{\rho}{2} \left\| \mathcal Z_1 - \mathcal X_0 + \mathcal U_4 \right\|_F^2 .
\end{equation*}

This yields
\begin{equation}
\mathcal Z_1 = \left[ \mathcal M \times_3 X_1^T + \rho (\mathcal X_0 - \mathcal U_4) \right] \times_3 \left( X_1^T X_1 + \rho I_k \right)^{-1}.
\label{eq:updateZ1}
\end{equation}

\item The seventh and eighth sub-problem read respectively

\begin{equation*}
\min_{\mathcal Z_2} \lambda_s \left\| \mathcal Z_2 \right\|_1 + \frac{\rho}{2} \left\| \mathcal Z_2 - \mathcal X_0 \times_1 D_m + \mathcal U_5 \right\|_F^2 ,
\end{equation*}

\begin{equation*}
\min_{\mathcal Z_3} \lambda_s \left\| \mathcal Z_3 \right\|_1 + \frac{\rho}{2} \left\| \mathcal Z_3 - \mathcal X_0 \times_2 D_n + \mathcal U_6 \right\|_F^2 .
\end{equation*}

Solutions for both problems are given by

\begin{equation*}
\mathcal Z_2 = \soft \left(\mathcal X_0 \times_1 D_m - \mathcal U_5, \frac{\lambda_s}{\rho} \right),
\label{eq:updateZ2}
\end{equation*}

\begin{equation*}
\mathcal Z_3 = \soft \left(\mathcal X_0 \times_2 D_n - \mathcal U_6, \frac{\lambda_s}{\rho} \right).
\label{eq:updateZ32}
\end{equation*}

\item The ninth sub-problem is another projection:

\begin{equation*}
\min_{\mathcal Z_4} \chi_+(\mathcal Z_4) + \frac{\rho}{2} \left\| \mathcal Z_4 - \mathcal X_0 + \mathcal U_7 \right\|_F^2.
\end{equation*}

The update is thus given by

\begin{equation*}
\mathcal Z_4 = \Pi_+ \left( \mathcal X_0 - \mathcal U_7 \right).
\label{eq:updateZ4}
\end{equation*}

\end{itemize}

The longest update is the one involving the tensor contraction, described in Equation \eqref{eq:updateZ1}. 
The computation of ${X_1}^T M^{(3)}$ indeed takes $O(k m n o)$ steps to compute.
The one involving the $\dct$ of the abundance matrix, described by Equation \eqref{eq:DCTupdateX0}, 
can be computed in $O (m n k \log(m n k))$ and is therefore usually lighter (depending on $o$).
This is also lighter that the DCT computations described in Section \ref{sec:TV}, where
Equation \eqref{eq:updateX3} required a $O (m n o \log(m n o))$ computation.

\section{Discussion and final remarks}

\begin{figure}[!ht]
\centering
\subfloat[\label{fig:nonoise}]{\includegraphics[width=.3\textwidth]{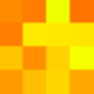}}\quad %
\subfloat[\label{fig:noisy}]{\includegraphics[width=.3\textwidth]{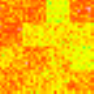}}\quad %
\subfloat[\label{fig:Mediandenoise}]{\includegraphics[width=.3\textwidth]{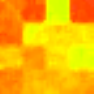}}\quad \newline%
\subfloat[\label{fig:Wienerdenoise}]{\includegraphics[width=.3\textwidth]{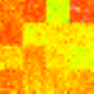}}\quad %
\subfloat[\label{fig:LeeSeungDenoise}]{\includegraphics[width=.3\textwidth]{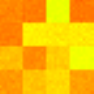}}\quad %
\subfloat[\label{fig:SPAdenoise}]{\includegraphics[width=.3\textwidth]{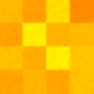}}\quad \newline%
\subfloat[\label{fig:TVdenoise0050}]{\includegraphics[width=.3\textwidth]{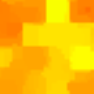}}\quad %
\subfloat[\label{fig:TVdenoise005001}]{\includegraphics[width=.3\textwidth]{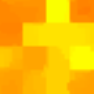}}\quad %
\subfloat[\label{fig:NMFTVdenoise201}]{\includegraphics[width=.3\textwidth]{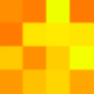}}\quad \newline%
\caption{Simulation results for band $132$. 
(\ref{fig:nonoise}) Simulated hypercube before adding noise,
(\ref{fig:noisy}) Simulated hypercube after adding noise,
(\ref{fig:Mediandenoise}) Median denoising,
(\ref{fig:Wienerdenoise}) Wiener denoising,
(\ref{fig:LeeSeungDenoise}) Lee and Seung's NMF used for denoising,
(\ref{fig:SPAdenoise}) SPA used for denoising,
(\ref{fig:TVdenoise0050}) TV denoising with $\lambda_s = 0.05$, and $\lambda_t = 0$,
(\ref{fig:TVdenoise005001}) TV denoising with $\lambda_s = 0.05$, and $\lambda_t = 0.01$,
(\ref{fig:NMFTVdenoise201}) NMF-TV denoising with $\lambda_s = 2.0$, and $\lambda_t = 0.1$.
}
\label{fig:simu}
\end{figure}

We were able to add TV regularization to the ADMM algorithm for NMF without slowing it down significantly.
Resulting unmixings seem to benefit from introduced regularization.

To illustrate the denoising techniques proposed in this article, we generate a synthetic hypercube, and apply the 
$3$-dimensional TV denoising introduced in Section \ref{sec:TV}, and the NMF-TV algorithm from Section \ref{sec:NMF}.
NMF was not designed to remove noise. Its goal is to recover non-negative factors from a matrix. However, these factors are usually
not unique, and can be quite different from the ones used to generate an example. Instead of assessing factors, we
therefore focus on their product. 

We generate a cube consisting of $16$ homogenous regions, each consisting of a randomly chosen mixture
of $5$ randomly chosen spectrums in the USGS 1995 
Library\footnote{Available online at \href{https://speclab.cr.usgs.gov/spectral.lib06/}{https://speclab.cr.usgs.gov/spectral.lib06/}}.
We then add white noise, and try out our techniques.

More precisely, we start by building true endmembers matrix $W_0 \in \R^{224 \times 5}$ by stacking up $5$
endmembers randomly selected amongst the $498$ ones available in the library.
We then build a first abundance matrix ${H_0} \in \R^{5 \times 16}$ where every column is randomly drawn from a 
Dirichlet distribution. The matrix product $W_0 H_0$ defines a matrix $M_0 \in \R^{ 224 \times 16 }$ whom we fold into
hypercube $\mathcal{M}_0\in \R^{ 4 \times 4 \times 224 }$.
To induce spatial regularity, we then replace each pixel $x$ of $\mathcal{M}_0$ by a $9 \times 9$ block of 
copies of $x$. We get a cube of size ${ 36 \times 36 \times 224 }$. We finally obtain hypercube $\mathcal{M}$ by adding 
gaussian noise.

\begin{equation*}
\mathcal{M} =  \mathcal{M}_0 +  \mathcal{N},
\end{equation*}

where every coefficient of $ \mathcal{N} \in \R^{ 36 \times 36 \times 224 }$ was randomly chosen according
to the centered normal distribution.

We compare proposed techniques with $4$ relatively common algorithms, a median filter, a Wiener filter, 
and two $NMF$ algorithms : the seminal multiplicative update scheme from Lee and Seung \cite{lee2001algorithms}, 
and the Successive Projection Algorithm (SPA) described by Gillis in \cite{gillis2014and}. The median filter 
replaces each pixel by the median of its $3 \times 3 \times 3$ neighborhood. The Wiener filter changes every 
pixel value $x$ according to
\begin{equation}
y= \begin{cases}
\frac{\sigma^2}{\sigma_x^2}m_x + \left( 1 - \frac{\sigma^2}{\sigma_x^2}\right) x \text{ if } \sigma_x^2 \geq \sigma^2, \\
m_x \text{ if } \sigma_x^2 < \sigma^2, \\
\end{cases}
\label{eq:wiener}
\end{equation}
where $m_x$ and $\sigma_x$ are the local mean and variance of the $3 \times 3 \times 3$ neighborhood of pixel $x$.

The median filter, Wiener filter and TV filter all give a result that we can directly compare to our synthesized example.
The NMF algorithms results consist of two factors $W$ and $H$. We use their product $W H$ for comparisons.

ADMM parameter $\rho$ can theoretically be set to any value without compromising convergence. In practice, 
a good choice for $\rho$ is important, as it has an influence on the algorithm's speed. We 
experimentally determined that setting $\rho = 10$ allows our algorithm to converge reasonably fast. 
Parameters $\lambda$ and $\rho$ controlling the desired strength of TV-regularization were also set 
experimentally.

Results are displayed in Figure \ref{fig:simu}. Only one slice of the resulting hypercubes is displayed. 
The median and Wiener filter are not able to recover precisely the boundaries of every region as can be seen
in Figure \ref{fig:Mediandenoise} and \ref{fig:Wienerdenoise}.
The classical NMF algorithms are not designed specifically with denoising in mind. However, their result displayed
in Figure \ref{fig:LeeSeungDenoise} and \ref{fig:SPAdenoise} show that they are quite successful at recovering sharp
boundaries between every region. They however fail at producing smooth regions.
Figure \ref{fig:TVdenoise0050} results from a spatial only TV regularization, whereas Figure \ref{fig:TVdenoise005001}
is regularized by spatial--spectral TV. A careful comparison of both figures reveals that adding a spectral component
to the TV denoising algorithm described in Section \ref{sec:TV} seems to slightly improve the result.
The NMF-TV algorithm recovers almost perfectly the original cube, as can be seen in Figure \ref{fig:NMFTVdenoise201}.

To conclude, we state an interesting direction for future research. The base ingredient we used throughout our paper
is a simple ADMM algorithm. We could instead try to use the Nesterov accelerated scheme described by 
Goldstein et al. in \cite{goldstein2014fast}. This would most likely speed up convergence.

\bibliographystyle{amsplain}
\bibliography{biblio}

\section*{Appendix} \label{app}

For an order $3$ tensor $\mathcal T \in \R^{m \times n \times o}$ the mode $1$-fiber corresponding to indices $(j,k)$ is the vector $(T_{i,j,k} )_{1\leq i \leq m}$, the $2$-fiber corresponding to $(i, k)$ is defined as $(T_{i,j,k} )_{1\leq j\leq n}$ and mode $3$-fibers are defined similarly as $(T_{i,j,k} )_{1\leq k\leq n}$. Mode $n$ matricization, also known as unfolding, of tensor $\mathcal T$ transforms a tensor of any order in a matrix, by stacking its $n$-fibers in a precise order. The mode $k$ unfolding of tensor $\mathcal T$ is denoted $T^{(k)}$. A mode $n$ tensor contraction with a matrix is the generalization of the usual matrix product. The notion is probably better understood through an example. Let $\mathcal T$ be an order $3$ tensor, and $X$ be a matrix. The mode $2$ contraction of $\mathcal T$ and $X$, denoted $\mathcal T \times_2 X$ is defined as

\begin{equation*}
\left( \mathcal T \times_2 X \right)_{i,l,k} = \sum_l T_{i,l,k} X_{l,j} .
\label{eq:contrex}
\end{equation*}

Tensor contractions, unfoldings, and matrix product are linked by the following property :

\begin{equation*}
\mathcal U = \mathcal T \times_k X \Leftrightarrow U^{(k)} = X T^{(k)}.
\label{eq:contrex2}
\end{equation*}

Multiple tensor contractions can be expressed with Kronecker products, thanks to one of its fundamental properties

\begin{equation*}
\left( B^T \otimes A \right) \mbox{vec}(X) = \mbox{vec}(A X B) = \mbox{vec}(C) .
\end{equation*}

This fact can be used to prove that the following propositions are equivalent :
\begin{enumerate}[(i)]
\item $\mathcal Y = \mathcal X \times_1 A^1 \times_2 A^2 \times_3 \cdots \times_n A^n,$
\item $Y^{(k)} = A^k X^{(k)} \left( A^n \otimes \cdots \otimes A^{k+1} \otimes A^{k-1} \otimes \cdots \otimes A^1 \right)^T ,$
\item $\mbox{vec}(Y^{(k)}) = \left( A^n \otimes \cdots \otimes A^{k+1} \otimes A^{k-1} \otimes \cdots \otimes A^1 \otimes A^k \right) \mbox{vec}(X^{(k)}).$
\end{enumerate}

Last property lets us derive successively to the vectorization of Equation \eqref{eq:Neumann3} in Section \ref{sec:TV} :

\begin{align*}
B^{(1)} = (I_{mno} + \rho L_o \otimes I_n \otimes I_m + \rho I_o \otimes L_n \otimes I_m + \rho I_o \otimes I_n \otimes L_m) X^{(1)}.
\end{align*}

Kronecker sums are usually defined by

\begin{equation*}
A \oplus B = A \otimes I_n + I_m \otimes B. 
\end{equation*}

This definition can be extended to

\begin{equation*}
A \oplus B \oplus C = A \otimes I_n \otimes I_k + I_m \otimes B \otimes I_k + I_m \otimes I_n \otimes C. 
\end{equation*}

Note that Kronecker sums are associative but do not commute.

A useful fact about Kronecker sums and products is that they behave well with eigenvalue decompositions.
Specifically, if $ A = Q_A D_A Q_A^T $ and $ B = Q_B D_B Q_B^T $ are the eigenvalue decompositions of matrices $A$ and $B$, then

\begin{align*}
& A \otimes B = (Q_A \otimes Q_B) (D_A \otimes D_B) (Q_A \otimes Q_B)^T, \\
& A \oplus B = (Q_A \otimes Q_B) (D_A \oplus D_B) (Q_A \otimes Q_B)^T,
\end{align*}

are the respective eigen-decompositions of $A \otimes B$ and $A \oplus B$.

\end{document}